\documentclass[12pt]{amsart}

\usepackage{amscd,latexsym,amsthm,amsfonts,amssymb,amsmath,amsxtra}
\usepackage[mathscr]{eucal}
\usepackage
[
colorlinks=true, 
linkcolor=blue,
urlcolor=blue,
bookmarks=true,
bookmarksopen=true,
citecolor=blue,
]
{hyperref}

\newcommand{\BA}{{\mathbb {A}}}

\newcommand{\BC}{{\mathbb {C}}}

\newcommand{\CA}{{\mathcal {A}}}

\newcommand{\CC}{{\mathcal {C}}}

\newcommand{\CE}{{\mathcal {E}}}
\newcommand{\CF}{{\mathcal {F}}}

\newcommand{\CM}{{\mathcal {M}}}

\newcommand{\CP}{{\mathcal {P}}}

\newcommand{\CR}{{\mathcal {R}}}
\newcommand{\CS}{{\mathcal {S}}}

\newcommand{\bs}{\backslash}

\newcommand{\GL}{{\mathrm{GL}}}

\newcommand{\Mat}{{\mathrm{Mat}}}

\newcommand{\PM}{{\mathrm{PM}}}

\newcommand{\Sp}{{\mathrm{Sp}}}

\def\tilsp{\mathop{\widetilde{\rm Sp}}\nolimits}

\newcommand{\wt}{\widetilde}

\def\tilb{{\widetilde{b}}}

\def\tilg{{\widetilde g}}

\def\tilm{{\widetilde m}}

\def\tilphi{{\widetilde \phi}}

\def\op{{\overline{p}}}
\def\oz{{\overline{z}}}

\def\sig{{\sigma}}

\def\tet{{\theta}}

\def\tiltet{{\widetilde{\tet}}}
\def\tiltau{{\widetilde{\tau}}}

\theoremstyle{plain}
\newtheorem{thm}{Theorem}[section]

\newtheorem{lem}[thm]{Lemma}
\newtheorem{prop}[thm]{Proposition}

\begin{document}
%------------------------------------------------------
\renewcommand{\theequation}{\arabic{equation}}
\numberwithin{equation}{section}

\title{Erratum to ``On the Non-vanishing of the Central Value of the Rankin-Selberg L-functions"}

\author{David Ginzburg}
\address{School of Mathematical Sciences\\
Sackler Faculty of Exact Sciences\\
Tel Aviv University\\
Ramat-Aviv, 69978 Israel}
\email{ginzburg@post.tau.ac.il}

\author{Dihua Jiang}
\address{School of Mathematics\\
University of Minnesota\\
Minneapolis, MN 55455, USA}
\email{dhjiang@math.umn.edu}

\author{Baiying Liu}
\address{Department of Mathematics\\
Purdue University\\
West Lafayette, IN 47906, USA}
\email{liu2053@purdue.edu}

\author{Stephen Rallis}
\address{Department of Mathematics\\
The Ohio State University\\
Columbus, OH 43210, USA}
\email{haar@math.ohio-state.edu}

%\dedicatory{Dedicated to Ilya I. Piatetski-Shapiro with admiration on the occasion of his 75th birthday}



\date{\today}
\begin{abstract}
    We complete the proof of Proposition 5.3 of \cite{GJR04}. 
\end{abstract}

\thanks{The second named author is partially supported by NSF grants DMS-1600685 and DMS-1901802. 
The third named author is partially supported by NSF grants DMS-1702218, DMS-1848058, and by start-up funds from the Department of Mathematics at Purdue University.}

%\tableofcontents

%\newpage

\maketitle

In this note, we complete the proof of Proposition 5.3 of \cite{GJR04} which is stated as follows. 

\begin{prop}[Proposition 5.3, \cite{GJR04}]\label{prop1}
If the period, defined in (5.2) of \cite{GJR04}, 
$$
{\CP}_{r,r-l}(\phi_\sig,\tilphi_\tiltau,\varphi_l),
$$
does not vanish for some given $\phi_\sig\in V_\sig$ and
$\tilphi_\tiltau\in V_\tiltau$, then the integral
\begin{equation*}
   \int_{K\times\Mat_{r-l,2r}(\BA)\times\PM_{2r,l}}
\phi_{\pi_\psi(\tiltau)\otimes\sig}(mk)
{\CF}^\psi(\tilphi_{\pi_\psi(\tiltau)\otimes\tiltau})
(mv^-(p_1)wk)dmdp_1dk, 
\end{equation*}
does not vanish for some choice of data 
$\phi_{\pi_\psi(\tiltau)\otimes\sig} \in \CA_{P_{2r,l}, \pi_{\psi}(\wt{\tau}) \otimes \sigma}$ and 
$\wt{\phi}_{\pi_\psi(\tiltau)\otimes\wt{\tau}} \in \wt{\CA}_{\wt{P}_{2r,r}, \pi_{\psi}(\wt{\tau}) \otimes \wt{\tau}, 0}$.
\end{prop}

The proof of Proposition 5.3 of \cite{GJR04} is reduced to the proof of the non-vanishing of the following integral (see \eqref{sec0equ4} below)
\begin{equation*}
\int_{\Mat_{r-l,2r}(\BA)\times\PM_{2r,l}}
\phi_{\pi_\psi(\tiltau)\otimes\sig}(m)
{\CF}^\psi(\tilphi_{\pi_\psi(\tiltau)\otimes\tiltau})
(mv^-(p_1)w)dmdp_1
\end{equation*}
for a proper set of sections, which was not complete in \cite{GJR04}. In this note, we complete this proof by proving Proposition \ref{prop2} below. Notation in the above proposition will be explained in Section \ref{notation}.

\subsection{Notation and the main result in Section 4 of \cite{GJR04}}\label{notation}

Let $F$ be a number field. 
We will use the notation from \cite{GJR04} freely. 
Let $\CA_{P_{2r,l}, \pi_{\psi}(\wt{\tau}) \otimes \sigma}$ be the set of functions
$$\phi: M_{2r,l}(F)U_{2r,l}(\BA) \bs \Sp_{4r+2l}(\BA) \rightarrow \BC,$$
such that $\phi$ is right $K_{\Sp_{4r+2l}}(\BA)$-finite, and for each $k \in K_{\Sp_{4r+2l}}(\BA)$, the function $\phi_k: m \mapsto \phi(mk)$, $m \in M_{2r,l}(\BA)$, belongs to $\pi_{\psi}(\wt{\tau}) \otimes \sigma$. 
For $\phi \in \CA_{P_{2r,l}, \pi_{\psi}(\wt{\tau}) \otimes \sigma}$, let 
$$\Phi(\cdot, s, \phi) = \phi(\cdot)\exp\langle s + \rho_{P_{2r,l}}, H_{P_{2r,l}}(\cdot)\rangle.$$
Then 
$$\left\{\Phi(\cdot, s, \phi): \phi \in \CA_{P_{2r,l}, \pi_{\psi}(\wt{\tau}) \otimes \sigma}\right\}$$
is equivalent to $I(s, \pi_{\psi}(\wt{\tau}) \otimes \sigma)$.
Similarly, let 
$\wt{\CA}_{\wt{P}_{2r,r}, \pi_{\psi}(\wt{\tau}) \otimes \wt{\tau}}$ be the set of functions
$$\wt{\phi}: \wt{M}_{2r,r}(F)U_{2r,r}(\BA) \bs \wt{\Sp}_{6r}(\BA) \rightarrow \BC,$$
such that $\wt{\phi}$ is right $K_{{\Sp}_{6r}}(\BA)$-finite, and for each $k \in K_{{\Sp}_{6r}}(\BA)$, the function $\wt{\phi}_k: \wt{m} \mapsto \wt{\phi}(\wt{m}k)$, $\wt{m} \in \wt{M}_{2r,r}(\BA)$, belongs to $\pi_{\psi}(\wt{\tau}) \otimes \wt{\tau}$. 
For $\wt{\phi} \in \wt{\CA}_{\wt{P}_{2r,r}, \pi_{\psi}(\wt{\tau}) \otimes \wt{\tau}}$, let 
$$\wt{\Phi}(\cdot, s, \wt{\phi}) = \wt{\phi}(\cdot)\gamma_{\psi}(\det(\cdot))
\exp\langle s + \rho_{\wt{P}_{2r,r}}, H_{\wt{P}_{2r,r}}(\cdot)\rangle.$$
Then 
$$\left\{\wt{\Phi}(\cdot, s, \wt{\phi}): \wt{\phi} \in \wt{\CA}_{\wt{P}_{2r,r}, \pi_{\psi}(\wt{\tau}) \otimes \wt{\tau}}\right\}$$
is equivalent to $\wt{I}(s, \pi_{\psi}(\wt{\tau}) \otimes \wt{\tau})$.

The goal of Section 4 of \cite{GJR04} is to compute the period
\begin{equation}\label{sec2equ9}
    \CP_{3r,r-l}(E_{\frac{1}{2}}(\cdot, \phi), \wt{E}_1(\cdot, \wt{\phi}), \varphi_{2r+l}) = 
    \int_{[\Sp_{4r+2l}]}E_{\frac{1}{2}}(g, \phi)\CF^{\psi}_{\varphi_{2r+l}}(\wt{E}_1(\cdot, \wt{\phi}))(g)dg,
\end{equation}
where, $\phi = \phi_{\pi_{\psi}(\wt{\tau}) \otimes \sigma} \in \CA_{P_{2r,l}, \pi_{\psi}(\wt{\tau}) \otimes \sigma}$, $E_{\frac{1}{2}}(\cdot, \phi)$ is the residue at $s=\frac{1}{2}$ of the following Eisenstein series
$$E(g, s, \phi) = \sum_{\gamma \in P_{2r,l}(F) \bs \Sp_{4r+2l}(F)} \Phi(\gamma g, s, \phi), \, g \in \Sp_{4r+2l}(\BA);$$
$\wt{\phi} =\wt{\phi}_{\pi_{\psi}(\wt{\tau}) \otimes \wt{\tau}} \in \wt{\CA}_{\wt{P}_{2r,r}, \pi_{\psi}(\wt{\tau}) \otimes \wt{\tau}}$, $\wt{E}_{1}(\cdot, \wt{\phi})$ is the residue at $s=1$ of the following Eisenstein series
$$\wt{E}(g, s, \wt{\phi}) = \sum_{\gamma \in P_{2r,r}(F) \bs \Sp_{6r}(F)} \wt{\Phi}(\gamma \wt{g}, s, \wt{\phi}), \, \wt{g} \in \wt{\Sp}_{6r}(\BA);$$
and 
$$\CF^{\psi}_{\varphi_{2r+l}}(\wt{E}_1(\cdot, \wt{\phi}))(g) = \int_{[V_{2r+l}]} \wt{\theta}^{\psi^{-1}}_{\varphi_{2r+l}}(\ell_{2r+l}(v)\wt{g})\wt{E}_1(v\wt{g}, \wt{\phi}) \psi_{r-l}(v)dv.$$
It turns out (see (4.8) of \cite{GJR04}) that the period \eqref{sec2equ9} is the residue at $s=\frac{1}{2}$ of the following period
\begin{equation}\label{sec2equ1}
    \CP_{3r,r-l}(\CE_1, \wt{E}_1(\cdot, \wt{\phi}), \varphi_{2r+l}) = 
    \int_{[\Sp_{4r+2l}]}\CE_1(g)\CF^{\psi}_{\varphi_{2r+l}}(\wt{E}_1(\cdot, \wt{\phi}))(g)dg,
\end{equation}
where, 
$$\CE_1(g) = \sum_{\gamma \in P_{2r,l}(F) \bs \Sp_{4r+2l}(F)} \Phi(\gamma g, s, \phi)(1-\tau_c(H(\gamma g))).$$
Recall from (4.4) of \cite{GJR04} that for $g = um(g)k \in \Sp_{4r+2l}(\BA)$ with $u \in U_{2r,l}(\BA)$, $m(g) \in M_{2r,l}(\BA)$ and $k\in K_{\Sp_{4r+2l}}(\BA)$, $H(g)=\lvert\det(m(g))\rvert$. 
We remark that all $\wt{\theta}$ occurred in Sections 4 and 5 of \cite{GJR04}, namely for the case of $r \geq l$, should be with respect to the character $\psi^{-1}$. 

Let $\Phi^c(\gamma g, s, \phi)=\Phi(\gamma g, s, \phi)(1-\tau_c(H(\gamma g)))$.
By \cite[Proposition 4.3]{GJR04}, the period \eqref{sec2equ1} is equal to 
\begin{equation}\label{sec2equ2}
    \int_{M(F)U(\BA)\bs\Sp_{4r+2l}(\BA)}
    \Phi^c(\gamma g, s, \phi)
    \int_{\Mat_{r-l,2r}(\BA)}
    \int_{[V_{r,l}]} \sum_{\xi_l \in F^l} \omega_{\psi^{-1}}(\ell(\overline{p}_2, \overline{z})\ell(\overline{p}_1)\wt{g})
\end{equation}
$$\varphi_{2r+l}((0,\xi_l)) \wt{E}_{1,\wt{P}_{2r,r}}(vv^-(p_1)w\wt{g}, \wt{\phi})\psi_{r-l}(v)dvdp_1dg,$$
where $\wt{E}_{1,\wt{P}_{2r,r}}(\wt{g}, \wt{\phi})$ is the constant term of the residue $\wt{E}_{1}(\wt{g}, \wt{\phi})$ along the maximal parabolic subgroup $\wt{P}_{2r,r}$, which equals $\wt{\CM}_1(\wt{\Phi})(\wt{g})$, the residue at $s=1$ of the intertwining operator 
$\wt{\CM}(w_{2r,r},s)(\wt{\Phi})(\wt{g})$ defined in Section 3.2 of \cite{GJR04}.
$\wt{\CM}(w_{2r,r},s)$ maps sections in the induced representation $\wt{I}(s, \pi_{\psi}(\wt{\tau}) \otimes \wt{\tau})$ to those in the induced representation 
$\wt{I}(-s, w_{2r,r}(\pi_{\psi}(\wt{\tau}) \otimes \wt{\tau}))$. Note that 
$\wt{\CM}_1(\wt{\Phi})(\wt{g})$ is not identically zero, and 
$w_{2r,r}(\pi_{\psi}(\wt{\tau}) \otimes \wt{\tau}) = \pi_{\psi}(\wt{\tau}) \otimes \wt{\tau}$ since $\pi_{\psi}(\wt{\tau})$ is self-dual. 

After applying the Iwasawa decomposition $\Sp_{4r+2l}(\BA)=P_{2r,l}(\BA)K$, $K=K_{\Sp_{4r+2l}}(\BA)$, we obtain the integral (4.31) of \cite{GJR04}, in which $\wt{\CM}_1(\wt{\Phi})$ belongs to $\wt{I}(-1, \pi_{\psi}(\wt{\tau}) \otimes \wt{\tau})$. 
Since the induced representation 
$\wt{I}(-1, \pi_{\psi}(\wt{\tau}) \otimes \wt{\tau})$ is 
reducible, the image $\wt{\CM}_1(\wt{\Phi})$ belongs to a proper 
subrepresentation of $\wt{I}(-1, \pi_{\psi}(\wt{\tau}) \otimes \wt{\tau})$. 
This is 
the key point that makes the original argument in the proof in \cite{GJR04} for Proposition 5.3 insufficient.  Denote the subrepresentation of $\wt{I}(-1, \pi_{\psi}(\wt{\tau}) \otimes \wt{\tau})$ consisting of the images of $\wt{\CM}_1$ by $\wt{I}_0(-1, \pi_{\psi}(\wt{\tau}) \otimes \wt{\tau})$.
Recall that $\wt{I}(-1, \pi_{\psi}(\wt{\tau}) \otimes \wt{\tau})$ is equivalent to
$$\left\{\wt{\Phi}(\cdot, -1, \wt{\phi}): \wt{\phi} \in \wt{\CA}_{\wt{P}_{2r,r}, \pi_{\psi}(\wt{\tau}) \otimes \wt{\tau}}\right\}.$$
Denote the subspace of $\wt{\CA}_{\wt{P}_{2r,r}, \pi_{\psi}(\wt{\tau}) \otimes \wt{\tau}}$ corresponding to $\wt{I}_0(-1, \pi_{\psi}(\wt{\tau}) \otimes \wt{\tau})$ by $\wt{\CA}_{\wt{P}_{2r,r}, \pi_{\psi}(\wt{\tau}) \otimes \wt{\tau}, 0}$.

In order to complete the proof of Proposition 5.3 in \cite{GJR04}, we find a technically more involved argument, which is not 
sensitive to which section to be taken in 
the subrepresentation $\wt{I}_0(-1, \pi_{\psi}(\wt{\tau}) \otimes \wt{\tau})$ or even in the whole induced representation $\wt{I}(-1, \pi_{\psi}(\wt{\tau}) \otimes \wt{\tau})$.

Since $\wt{\CM}_1(\wt{\Phi}) \in \wt{I}_0(-1, \pi_{\psi}(\wt{\tau}) \otimes \wt{\tau})$, 
from the discussion above, there exists
$\tilphi_{\pi_\psi(\tiltau)\otimes\tiltau} \in \wt{\CA}_{\wt{P}_{2r,r}, \pi_{\psi}(\wt{\tau}) \otimes \wt{\tau}, 0}$, such that 
$$\wt{\CM}_1(\wt{\Phi})=\tilphi_{\pi_\psi(\tiltau)\otimes\tiltau}\,\gamma_{\psi}(\det)
\exp\langle -1 + \rho_{\wt{P}_{2r,r}}, H_{\wt{P}_{2r,r}}\rangle.$$ 
It follows that 
\begin{eqnarray*}
&&\wt{\CM}_1(\wt{\Phi})(\wt{m}(a,b)vv^-(p_1)wk)\\
=&&
\lvert \det a \rvert^{-1+2r+\frac{1}{2}} \gamma_{\psi}(\det a)
\tilphi_{\pi_\psi(\tiltau)\otimes\tiltau}(\wt{m}(a,b)vv^-(p_1)wk).
\end{eqnarray*}
After carrying out the calculations from (4.32) to (4.35) of \cite{GJR04}, one obtains in Theorem 4.4 of \cite{GJR04} that the period \eqref{sec2equ9} is equal to a product that a constant $\frak{c}$ times the integral (4.35) of \cite{GJR04} which is given by 
\begin{equation}\label{sec2equ8}
   \int_{K\times\Mat_{r-l,2r}(\BA)\times\PM_{2r,l}}
\phi_{\pi_\psi(\tiltau)\otimes\sig}(mk)
{\CF}^\psi(\tilphi_{\pi_\psi(\tiltau)\otimes\tiltau})
(mv^-(p_1)wk)dmdp_1dk, 
\end{equation}
where the function ${\CF}^\psi(\tilphi_{\pi_\psi(\tiltau)\otimes\tiltau})(mv^-(p_1)wk)$ is defined
as in (4.33) of \cite{GJR04} by 
\begin{equation}\label{sec0equ1}
\int_{[V_{r,l}]}\tiltet_{\varphi_{2r+l},l}^{\psi^{-1}}(\ell(\op_2,\oz)\ell(\op_1)\tilm(b)k)
\tilphi_{\pi_\psi(\tiltau)\otimes\tiltau}(v\tilm(a,b)v^-(p_1)wk)\psi_{r-l}(v)dv,
\end{equation}
the integration domain $\PM_{2r,l}$ is given by 
\begin{equation}\label{sec2equ6}
    \PM_{2r,l}:=(Z_{\GL_{2r}}(\BA)\GL_{2r}(F)\bs\GL_{2r}(\BA))\times
(\Sp_{2l}(F)\bs\Sp_{2l}(\BA))
\end{equation}
as in (4.34) of \cite{GJR04},
and
\begin{equation}\label{sec2equ5}
    \frak{c} = \frac{\mathrm{vol}(\BA^1 \bs F^{\times})}{2rd}
\end{equation}
with $d$ being the number of the real archimedean places of the number field $F$.
Recall from 
(4.29) of \cite{GJR04} that $\tiltet_{\varphi_{2r+l},l}^{\psi^{-1}}$ is defined as follows:
\begin{equation}\label{sec2equ10}
    \tiltet_{\varphi_{2r+l},l}^{\psi^{-1}}(\ell(\bar{p}_2,\bar{z})\ell(\bar{p}_1)\wt{g}) = 
    \sum_{\xi_l \in F^l} \omega_{\psi^{-1}}(\ell(\bar{p}_2,\bar{z})\ell(\bar{p}_1)\wt{g})\varphi_{2r+l}((0,\xi_l)).
\end{equation}

Recall that for each $k \in K_{\Sp_{4r+2l}}(\BA)$,
$\phi \in \CA_{P_{2r,l}, \pi_{\psi}(\wt{\tau}) \otimes \sigma}$, $\phi_k: m \mapsto \phi(mk)$ belongs to $\pi_{\psi}(\wt{\tau}) \otimes \sigma$, and for each $k \in K_{\Sp_{6r}}(\BA)$, $\wt{\phi} \in \wt{\CA}_{\wt{P}_{2r,r}, \pi_{\psi}(\wt{\tau}) \otimes \wt{\tau}}$, $\wt{\phi}_k: \wt{m} \mapsto \wt{\phi}(\wt{m}k)$ belongs to $\pi_{\psi}(\wt{\tau}) \otimes \wt{\tau}$.
From now on, we consider special sections $\phi_{\pi_\psi(\tiltau)\otimes\sig} \in \CA_{P_{2r,l}, \pi_{\psi}(\wt{\tau}) \otimes \sigma}$ and 
$\wt{\phi}_{\pi_\psi(\tiltau)\otimes\wt{\tau}} \in \wt{\CA}_{\wt{P}_{2r,r}, \pi_{\psi}(\wt{\tau}) \otimes \wt{\tau}, 0}$ which are factorizable and have the property that 
$\phi_{\pi_\psi(\tiltau)\otimes\sig, 1}$ and 
$\wt{\phi}_{\pi_\psi(\tiltau)\otimes\wt{\tau}, w}$ are decomposable in $\pi_{\psi}(\wt{\tau}) \otimes \sigma$ and $\pi_{\psi}(\wt{\tau}) \otimes \wt{\tau}$, respectively:
\begin{equation}\label{sec2equ3}
    \phi_{\pi_\psi(\tiltau)\otimes\sig, 1} = \phi_{\pi_\psi(\tiltau)} \otimes \phi_{\sigma} \in \pi_{\psi}(\wt{\tau}) \otimes \sigma,
\end{equation}
\begin{equation}\label{sec2equ4}
    \wt{\phi}_{\pi_\psi(\tiltau)\otimes\wt{\tau}, w} = \overline{\phi}_{\pi_\psi(\tiltau)} \otimes \wt{\phi}_{\wt{\tau}} \in \pi_{\psi}(\wt{\tau}) \otimes \wt{\tau},
\end{equation}
where 
$\overline{\phi}_{\pi_\psi(\tiltau)}$ is the complex conjugate of 
${\phi}_{\pi_\psi(\tiltau)}$, and $w$ is the following Weyl element on Page 696 of \cite{GJR04},
$$w=\begin{pmatrix}
0 & I_{2r} & 0 & 0 & 0\\
I_{r-l} & 0 & 0 & 0 & 0\\
0 & 0 & I_{2l} & 0 & 0\\
0 & 0 & 0 & 0 & I_{r-l}\\
0 & 0 & 0 & I_{2r} & 0
\end{pmatrix}.$$
From the definitions of $\CA_{P_{2r,l}, \pi_{\psi}(\wt{\tau}) \otimes \sigma}$ and 
$\wt{\CA}_{\wt{P}_{2r,r}, \pi_{\psi}(\wt{\tau}) \otimes \wt{\tau}, 0}$, it is clear that these special sections are dense in $\CA_{P_{2r,l}, \pi_{\psi}(\wt{\tau}) \otimes \sigma}$ and 
$\wt{\CA}_{\wt{P}_{2r,r}, \pi_{\psi}(\wt{\tau}) \otimes \wt{\tau}, 0}$, respectively. 
Denote the sets of these special sections by 
$\CA^{ss}_{P_{2r,l}, \pi_{\psi}(\wt{\tau}) \otimes \sigma}$ and 
$\wt{\CA}^{ss}_{\wt{P}_{2r,r}, \pi_{\psi}(\wt{\tau}) \otimes \wt{\tau}, 0}$, respectively.

\subsection{Proof of Proposition 5.3 of \cite{GJR04}}\label{proof}

We repeat the proof for Proposition 5.3 in \cite{GJR04} and point out the place that needs a more technical argument, which is now taken care of 
by Proposition \ref{prop2} below.

Recall from (5.2), (2.11), and (2.13) of \cite{GJR04} that the period
${\CP}_{r,r-l}(\phi_\sig,\tilphi_\tiltau,\varphi_l)$ equals 
\begin{equation}\label{sec2equ11}
    \int_{[\Sp_{2l}]}\phi_\sig(g)\int_{[V_{r,l}]}\tiltet^{\psi^{-1}}_{\varphi_l}(\ell_l(v)\tilg)
\tilphi_\tiltau(v\tilg)\psi_{r-l}(v)dvdg.
\end{equation}
It defines a continuous functional on the space of
$$
V_\sig\otimes\left(\wt{\Theta}^{\psi^{-1}}_l\otimes V_\tiltau\right)^{V_{r,l}, \psi_{r-l}},
$$
where $\wt{\Theta}^{\psi^{-1}}_l$ is the space generated by the theta functions
$\wt{\theta}^{\psi^{-1}}_{\varphi_l}$ with $\varphi_l\in{\CS}(\BA^l)$ and
$(\wt{\Theta}^{\psi^{-1}}_l\otimes V_\tiltau)^{V_{r,l},\psi_{r-l}}$ is
the space generated by the Fourier-Jacobi coefficients of automorphic forms in
$\tiltau$.

It is clear that ${\CS}(\BA^{2r+l})={\CS}(\BA^{2r})\widehat{\otimes}{\CS}(\BA^l)$. If we
take $\varphi_{2r+l}=\varphi_{2r}\otimes\varphi_l$ (separation of variables),
then we have
$$
\tiltet_{\varphi_{2r+l},l}^{\psi^{-1}}(\ell(\op_2,\oz)\ell(\op_1)\tilg)
=
\varphi_{2r}(\ell(\op_1))\cdot\tiltet^{\psi^{-1}}_{\varphi_l}(\ell(\op_2,\oz)\tilg)
$$
for $\tilg\in\tilsp_{2l}(\BA)$ (see \eqref{sec2equ10} for the definition of $\tiltet_{\varphi_{2r+l},l}^{\psi^{-1}}$). 
For any fixed $\varphi_{2r}\in {\CS}(\BA^{2r})$,
we consider all Bruhat-Schwartz functions
$$
\varphi_{2r+l}=\varphi_{2r}\otimes\varphi_l\in{\CS}(\BA^{2r+l}),
$$
with $\varphi_l\in{\CS}(\BA^l)$. It follows that the space generated by
$\tiltet_{\varphi_{2r+l},l}^{\psi^{-1}}(\ell(\op_2,\oz)\tilg)$ (with a fixed
$\varphi_{2r}\in {\CS}(\BA^{2r})$ and all $\varphi_l\in{\CS}(\BA^l)$) is
the same as the space $\wt{\Theta}^{\psi^{-1}}_l$ (generated by all
$\tiltet^{\psi^{-1}}_{\varphi_l}(\ell(\op_2,\oz)\tilg)$) as automorphic representations of
the Jacobi group $\tilsp_{2l}(\BA)\ltimes H_l(\BA)$, where $H_l$ is the
Heisenberg group generated by all $\ell(\op_2,\oz)$. In the following we may assume that
$\varphi_{2r}$ is supported in a small neighborhood of zero.

It follows that the
non-vanishing of the period ${\CP}_{r,r-l}(\phi_\sig,\tilphi_\tiltau,\varphi_l)$
is equivalent to the non-vanishing of the following integral
$$
\int_{[\Sp_{2l}]}\phi_\sig(b)\int_{[V_{r,l}]}
\tiltet_{\varphi_{2r+l},l}^{\psi^{-1}}(\ell(\op_2,\oz)\tilb)\tilphi_\tiltau(v\tilb)
\psi_{r-l}(v)dvdb.
$$
On the other hand, it is clear that the integral
$$
\int_{Z_{\GL_{2r}}(\BA)\GL_{2r}(F)\bs\GL_{2r}(\BA)}
\phi_{\pi_\psi(\tiltau)}(a)
\overline{\phi_{\pi_\psi(\tiltau)}}(a)da
$$
is not zero for any choice of nonzero $\phi_{\pi_\psi(\tiltau)}$ and
$\overline{\phi_{\pi_\psi(\tiltau)}}$.
Hence, combining 
the above
two non-vanishing integrals, we obtain that the integral 
\begin{equation}\label{sec0equ2}
\int_{\PM_{2r,l}}
\phi_{\pi_\psi(\tiltau)\otimes\sig}(m)
{\CF}^\psi(\tilphi_{\pi_\psi(\tiltau)\otimes\tiltau})(mw)dm,
\end{equation}
does not vanish for some choice of data 
$\phi_{\pi_\psi(\tiltau)\otimes\sig} \in \CA^{ss}_{P_{2r,l}, \pi_{\psi}(\wt{\tau}) \otimes \sigma}$ and 
$\wt{\phi}_{\pi_\psi(\tiltau)\otimes\wt{\tau}} \in \wt{\CA}^{ss}_{\wt{P}_{2r,r}, \pi_{\psi}(\wt{\tau}) \otimes \wt{\tau}, 0}$, 
where $\PM_{2r,l}$ is as in \eqref{sec2equ6}, and for $m=m(a,b)$, ${\CF}^\psi(\tilphi_{\pi_\psi(\tiltau)\otimes\tiltau})(mw)$ is defined by
\begin{equation}\label{sec2equ7}
\int_{[V_{r,l}]}\tiltet_{\varphi_{2r+l},l}^{\psi^{-1}}(\ell(\op_2,\oz)\tilm(b))
\tilphi_{\pi_\psi(\tiltau)\otimes\tiltau}(\wt{m}(a)v\tilm(b)w)\psi_{r-l}(v)dv.
\end{equation}

We claim that 
for any choice of  $\phi_{\pi_\psi(\tiltau)\otimes\sig} \in \CA^{ss}_{P_{2r,l}, \pi_{\psi}(\wt{\tau}) \otimes \sigma}$, there exists $\wt{\phi}_{\pi_\psi(\tiltau)\otimes\wt{\tau}} \in \wt{\CA}^{ss}_{\wt{P}_{2r,r}, \pi_{\psi}(\wt{\tau}) \otimes \wt{\tau}, 0}$ such that the integral
\eqref{sec0equ2} does not vanish. 
Indeed, from
the discussion above and from the definitions of the special sections in 
$\CA^{ss}_{P_{2r,l}, \pi_{\psi}(\wt{\tau}) \otimes \sigma}$ and $\wt{\CA}^{ss}_{\wt{P}_{2r,r}, \pi_{\psi}(\wt{\tau}) \otimes \wt{\tau}, 0}$, 
it suffices to show that 
for any choice of $\phi_{\sigma} \in \sigma$, there exists $\wt{\phi}_{\wt{\tau}} \in \wt{\tau}$ such that the period  ${\CP}_{r,r-l}(\phi_\sig,\tilphi_\tiltau,\varphi_l)$ is nonzero. This follows from the fact that since $\sigma$ is irreducible, if the period  ${\CP}_{r,r-l}(\phi_\sig,\tilphi_\tiltau,\varphi_l)$ is nonzero for some choice of 
$\phi_{\sigma} \in \sigma$ and $\wt{\phi}_{\wt{\tau}} \in \wt{\tau}$, then 
the whole $\sigma$ occurs in the descent module of $\wt{\tau}$ (for the definition of descent modules see \cite[Chapter 3]{GRS11}).

Next we consider the following inner integration from Proposition \ref{prop1}:
\begin{equation}\label{sec0equ4}
\int_{\Mat_{r-l,2r}(\BA)\times\PM_{2r,l}}
\phi_{\pi_\psi(\tiltau)\otimes\sig}(m)
{\CF}^\psi(\tilphi_{\pi_\psi(\tiltau)\otimes\tiltau})
(mv^-(p_1)w)dmdp_1.
\end{equation}
Recall from Page 697 of \cite{GJR04} that the element $v^-(p_1)$ belongs to a unipotent subgroup of $\Sp_{6r}$ consisting of elements of the form
$$
v^-(p_1)=\begin{pmatrix}
I_{2r}&&&&\\
p_1&I_{r-l}&&&\\
&&I_{2l}&&\\
&&&I_{r-l}&\\
&&&p_1^*&I_{2r}\end{pmatrix}.
$$
By Proposition \ref{prop2} below and the claim above, 
for any choice of  $\phi_{\pi_\psi(\tiltau)\otimes\sig} \in \CA^{ss}_{P_{2r,l}, \pi_{\psi}(\wt{\tau}) \otimes \sigma}$, there exists $\wt{\phi}_{\pi_\psi(\tiltau)\otimes\wt{\tau}} \in \wt{\CA}_{\wt{P}_{2r,r}, \pi_{\psi}(\wt{\tau}) \otimes \wt{\tau}, 0}$ such that the integral
\eqref{sec0equ4} does not vanish.
This is the place where the original 
argument in the proof of Proposition 5.3 of \cite{GJR04} is not complete. 
Proposition \ref{prop2} will be proved in Sections \ref{idea} -- \ref{general case}. 

In order to prove finally the integral
\begin{equation*}\label{sec0equ5}
\int_{K\times\Mat_{r-l,2r}(\BA)\times\PM_{2r,l}}
\phi_{\pi_\psi(\tiltau)\otimes\sig}(mk)
{\CF}^\psi(\tilphi_{\pi_\psi(\tiltau)\otimes\tiltau})
(mv^-(p_1)wk)dmdp_1dk
\end{equation*}
is nonzero for some choice of data, for $k \in K$, we set
$$
\Psi(k):=
\int_{\Mat_{r-l,2r}(\BA)\times\PM_{2r,l}}
\phi_{\pi_\psi(\tiltau)\otimes\sig}(mk)
{\CF}^\psi(\tilphi_{\pi_\psi(\tiltau)\otimes\tiltau})
(mv^-(p_1)wk)dmdp_1.
$$
According to the discussion above, under the assumption of the proposition, for any choice of  $\phi_{\pi_\psi(\tiltau)\otimes\sig} \in \CA^{ss}_{P_{2r,l}, \pi_{\psi}(\wt{\tau}) \otimes \sigma}$, there exists $\wt{\phi}_{\pi_\psi(\tiltau)\otimes\wt{\tau}} \in \wt{\CA}_{\wt{P}_{2r,r}, \pi_{\psi}(\wt{\tau}) \otimes \wt{\tau}, 0}$ such that $\Psi(k)$ is nonzero at $k=1$, the identity. 
Since the special sections in  $\CA^{ss}_{P_{2r,l}, \pi_{\psi}(\wt{\tau}) \otimes \sigma}$ are dense in 
$\CA_{P_{2r,l}, \pi_{\psi}(\wt{\tau}) \otimes \sigma}$, we have the freedom on the $K$-support of the factorizable section $\phi_{\pi_\psi(\tiltau)\otimes\sig}$.

Therefore, at non-archimedean ramified local places $v$, we can choose a small support $\Omega_v \subset K_v$ of
$\phi_{\pi_\psi(\tiltau)\otimes\sig}$ near the identity, such that
$$
\Psi(k_\infty\cdot k_v)=\Psi(k_\infty).
$$
At the archimedean local places $v$, by using the continuity at $k=1$ of $\Psi(k)$, there is a small support
$\Omega_\infty\subset K_\infty$ for
$\phi_{\pi_\psi(\tiltau)\otimes\sig}$ such that the integral
$$
\int_K\Psi(k)dk=c_f\cdot\int_{\Omega_\infty}\Psi(k_\infty)dk_\infty\neq 0,
$$
with a constant $c_f$ depending on the ramified finite local places. 

This completes the proof of Proposition 5.3 of \cite{GJR04}, up to proving Proposition \ref{prop2} below.
\qed

\begin{prop}\label{prop2}
For any choice of data  $\phi_{\pi_\psi(\tiltau)\otimes\sig} \in \CA^{ss}_{P_{2r,l}, \pi_{\psi}(\wt{\tau}) \otimes \sigma}$, there exists  $\wt{\phi}_{\pi_\psi(\tiltau)\otimes\wt{\tau}} \in \wt{\CA}_{\wt{P}_{2r,r}, \pi_{\psi}(\wt{\tau}) \otimes \wt{\tau}, 0}$ such that 
the integral \eqref{sec0equ4},
which equals
$$
\int_{\Mat_{r-l,2r}(\BA)\times\PM_{2r,l}}
\phi_{\pi_\psi(\tiltau)\otimes\sig}(m)
{\CF}^\psi(\tilphi_{\pi_\psi(\tiltau)\otimes\tiltau})
(mv^-(p_1)w)dmdp_1,
$$
does not vanish.
\end{prop}

The proof of this proposition will be given in following sections. 

\subsection{The idea for proving Proposition \ref{prop2}}\label{idea}

In this section, we briefly introduce the idea for proving Proposition \ref{prop2}. 
First, we recall a lemma from \cite{GRS11}, which plays the same role as \cite[Lemma 4.2]{GJR04}. 

Let $H$ be any $F$-quasisplit classical group, including the general linear group.
Let $C$ be an $F$-subgroup of a maximal unipotent subgroup of $H$, and let $\psi_C$ be a non-trivial character of $[C] = C(F) \bs C(\BA)$.
$X, Y$ are two unipotent $F$-subgroups, satisfying the following conditions:
\begin{itemize}
\item[(1)] $X$ and $Y$ normalize $C$;
\item[(2)] $X \cap C$ and $Y \cap C$ are normal in $X$ and $Y$, respectively, $(X \cap C) \bs X$ and $(Y \cap C) \bs Y$ are abelian;
\item[(3)] $X(\BA)$ and $Y(\BA)$ preserve $\psi_C$;
\item[(4)] $\psi_C$ is trivial on $(X \cap C)(\BA)$ and $(Y \cap C)(\BA)$;
\item[(5)] $[X, Y] \subset C$;
\item[(6)]  there is a non-degenerate pairing $(X \cap C)(\BA) \times (Y \cap C)(\BA) \rightarrow \BC^*$, given by $(x,y) \mapsto \psi_C([x,y])$, which is
multiplicative in each coordinate, and identifies $(Y \cap C)(F) \bs Y(F)$ with the dual of
$
X(F)(X \cap C)(\BA) \bs X(\BA),
$
and
$(X \cap C)(F) \bs X(F)$ with the dual of
$
Y(F)(Y \cap C)(\BA) \bs Y(\BA).
$
\end{itemize}

Let $B =CY$ and $D=CX$, and extend $\psi_C$ trivially to characters of $[B]=B(F)\bs B(\BA)$ and $[D]=D(F)\bs D(\BA)$,
which will be denoted by $\psi_B$ and $\psi_D$ respectively. When there is no confusion, we may denote $\psi_B$ and $\psi_D$ all by $\psi_C$. 

\begin{lem}[Lemma 7.1 of \cite{GRS11}]\label{lem1}
Assume that the quadruple $(C, \psi_C, X, Y)$ satisfies all the above conditions. Let $f$ be an automorphic form on $H(\BA)$. Then for any $g \in H(\BA)$, 
$$\int_{[B]} f(vg) \psi_B^{-1}(v) dv =
\int_{(Y\cap C)(\BA) \backslash Y(\BA)}\int_{[D]} f(uyg) \psi_D^{-1}(u) du dy.$$
The right hand side of the the above equality is convergent in the sense
$$
\int_{(Y\cap C)(\BA) \backslash Y(\BA)} \vert \int_{[D]} f(uyg) \psi_D^{-1}(u) du \vert dy < \infty,$$
and this convergence is uniform as $g$ varies in compact subsets of $H(\BA)$. 
\end{lem}

\textbf{The idea for proving Proposition \ref{prop2}.}

First, based on the discussion in Section \ref{proof}, for any choice of $\phi_{\pi_\psi(\tiltau)\otimes\sig} \in \CA^{ss}_{P_{2r,l}, \pi_{\psi}(\wt{\tau}) \otimes \sigma}$, there exists $\wt{\phi}_{\pi_\psi(\tiltau)\otimes\wt{\tau}} \in \wt{\CA}^{ss}_{\wt{P}_{2r,r}, \pi_{\psi}(\wt{\tau}) \otimes \wt{\tau}, 0}$ such that the integral \eqref{sec0equ2}:
\begin{equation}\label{sec0equ20}
\int_{\PM_{2r,l}}
\phi_{\pi_\psi(\tiltau)\otimes\sig}(m)
{\CF}^\psi(\tilphi_{\pi_\psi(\tiltau)\otimes\tiltau})(mw)dm,
\end{equation}
does not vanish, 
where $\PM_{2r,l}$ is as in \eqref{sec2equ6}, and for $m=m(a,b)$,
${\CF}^\psi(\tilphi_{\pi_\psi(\tiltau)\otimes\tiltau})(mw)$ is defined in \eqref{sec2equ7}.

The proof of Proposition \ref{prop2} briefly consists of the following 4 steps. 
\begin{enumerate}
    \item Reversing the calculations from (4.29) -- (4.35) and reversing the step of taking Fourier expansion of $\wt{E}_1$ along $[\CC_1]$ as in \cite[Section 4]{GJR04}, we can  transform the integral \eqref{sec0equ20} to a nonzero constant times the residue at $s=\frac{1}{2}$ of a multiple integral over $[M]$ and $[V^{(1)}]$ (see \eqref{prop3equ5}, \eqref{sec0equ9} below).
    
    \item Note that $V^{(1)}=\prod_{n=2}^{r-l}\CC_nV^{(0)}$, we consider the integral over $[V^{(1)}]$.
Applying Lemma \ref{lem1} repeatedly to exchange roots from $\prod_{i=2}^{r-l}\CC_i$ to 
$\prod_{j=1}^{r-l-1}\CR_j$, we obtain a multiple integral over 
$\prod_{i=2}^{r-l}\CC_i(\BA)$, 
$[\prod_{j=1}^{r-l-1}\CR_j]$ and 
    $[V^{(0)}]$ (see \eqref{prop3equ7}, \eqref{sec0equ11} below). Combining with the outer integral over $[M]$, after changing of variables, we obtain a non-vanishing multiple integral over 
    $\prod_{i=2}^{r-l}\CC_i(\BA)$, $[M]$, $[\prod_{j=1}^{r-l-1}\CR_j]$ and 
    $[V^{(0)}]$ (see \eqref{prop3equ9}, \eqref{sec0equ12} below).
    Then we drop the outer integral over $\prod_{i=2}^{r-l}\CC_i(\BA)$.  
    Clearly the inner multiple integral over $[M]$, $[\prod_{j=1}^{r-l-1}\CR_j]$ and 
    $[V^{(0)}]$ is non-vanishing. 
    
    \item Then we consider the non-vanishing inner multiple integral over $[\prod_{j=1}^{r-l-1}\CR_j]$ and $[V^{(0)}]$. 
Applying Lemma \ref{lem1} repeatedly to exchange roots 
from 
$\prod_{j=1}^{r-l-1}\CR_j$ to  $\prod_{i=2}^{r-l}\CC_i$, we obtain a multiple integral over
$\prod_{j=1}^{r-l-1}\CR_j(\BA)$,  $[\prod_{i=2}^{r-l}\CC_i]$ and 
    $[V^{(0)}]$ (see \eqref{prop3equ12}, \eqref{sec0equ15} below). 
    Combining with the outer integral over $[M]$, after changing of variables, we obtain a non-vanishing multiple integral over 
   $\prod_{j=1}^{r-l-1}\CR_j(\BA)$, $[M]$, $[\prod_{i=2}^{r-l}\CC_i]$ and 
    $[V^{(0)}]$  (see \eqref{prop3equ18}, \eqref{sec0equ17} below).
    Note that $\prod_{i=2}^{r-l}\CC_iV^{(0)}=V^{(1)}$.
    
    \item By choosing appropriate $\varphi_{2r} \in \mathcal{S}(\BA^{2r})$, we obtain a non-vanishing multiple integral over 
   $\prod_{j=1}^{r-l}\CR_j(\BA)$, $[M]$ and 
    $[V^{(1)}]$  (see \eqref{prop3equ17}, \eqref{sec0equ19} below).
    Note that $\prod_{j=1}^{r-l}\CR_j(\BA)$ is exactly 
    the group $\Mat_{r-l,2r}(\BA)$. After taking Fourier expansion of $\wt{E}_1$ along $[\CC_1]$ and the calculations from (4.29) -- (4.35) as in \cite[Section 4]{GJR04}, we obtain a non-vanishing integral which is exactly a product of a nonzero constant with the integral in \eqref{sec0equ4}, for $\phi_{\pi_\psi(\tiltau)\otimes\sig}$ and some right translation of  $\wt{\phi}_{\pi_\psi(\tiltau)\otimes\wt{\tau}}$.
\end{enumerate}

\subsection{Proof of Proposition \ref{prop2}: special case $l=1, r=3$}\label{special case}

In this section, we prove Proposition \ref{prop2} for the first non-trivial case: $l=1, r=3$. 

We start from the non-vanishing integral \eqref{sec0equ20}.
Reversing the calculations from (4.29) -- (4.35) as in \cite[Section 4]{GJR04}, the integral \eqref{sec0equ20} is equal to 
$\frac{1}{\frak{c}}$ times the residue at $s=\frac{1}{2}$ of 
\begin{equation}\label{prop3equ4}
\int_{[M]}\Phi^c(m,s,\phi) 
\lvert \det a \rvert^{-\frac{21}{2}}
    \int_{[V_{3,1}]}
    \sum_{\xi_1 \in F^1} \omega_{\psi^{-1}}(\ell(\bar{p}_2,\bar{z})\wt{m})
\end{equation}
\begin{equation*}
   \varphi_{7}((0,\xi_1))\wt{E}_{1,\wt{P}_{6,3}}(v\wt{m}w,\wt{\phi})\psi_{2}(v)dvdm,
\end{equation*}
where $\frak{c}$ is as in \eqref{sec2equ5}. 

Reversing the step of taking Fourier expansion of $\wt{E}_1$ along $[\CC_1]$ as in \cite[Section 4]{GJR04}, the integral \eqref{prop3equ4} is equal to 
\begin{equation}\label{prop3equ5}
\int_{[M]}\Phi^c(m,s,\phi)
\lvert \det a \rvert^{-\frac{21}{2}}
    \int_{[V^{(1)}]}
    \sum_{\xi_1 \in F^1} \omega_{\psi^{-1}}(\ell(\bar{p}_2,\bar{z})\wt{m})
\end{equation}
\begin{equation*}
   \varphi_{7}((0,\xi_1))\wt{E}_{1}(v^{(1)}\wt{m}w,\wt{\phi})
   \psi_{2}(v^{(1)})dv^{(1)}dm.
\end{equation*}
Recall that $V^{(1)}$ consists of elements of the type
$$v^{(1)}=\begin{pmatrix}
I_6 & q & y & p_3^* & z'\\
    & n & p_2 & z & p_3\\
    &   & I_2 & p_2^* & y^*\\
    &   &     & n^* & q^*\\
    &   &     &     & I_6
\end{pmatrix},$$
where $q \in \Mat_{6,2}$ with the first column being zero. 

Recall that $V^{(0)}$ consists of elements of the type
$$v^{(0)}=\begin{pmatrix}
I_6 & 0 & y & p_3^* & z'\\
    & n & p_2 & z & p_3\\
    &   & I_2 & p_2^* & y^*\\
    &   &     & n^* & 0\\
    &   &     &     & I_6
\end{pmatrix}.$$
And for $1 \leq t \leq 2$,
$$\CC_t = \left\{
\begin{pmatrix}
I_6 & q & 0 & 0 & 0\\
    & I_2 & 0 & 0 & 0\\
    &   & I_2 & 0 & 0\\
    &   &     & I_2 & q^*\\
    &   &     &     & I_6
\end{pmatrix}: q \in \Mat_{6,2}, q_{i,j}=0, j \neq t \right\},$$
$$\CR_t = \left\{
\begin{pmatrix}
I_6 & 0 & 0 & 0 & 0\\
  p  & I_2 & 0 & 0 & 0\\
    &   & I_2 & 0 & 0\\
    &   &     & I_2 & 0\\
    &   &     &   p^*  & I_6
\end{pmatrix}: p \in \Mat_{2,6}, q_{i,j}=0, i \neq t \right\}.$$

Note that $V^{(1)}=\CC_2V^{(0)}$. Next, we consider the integral over $[V^{(1)}]$ and apply Lemma \ref{lem1} to exchange the roots from $\CC_2$ to 
$\CR_1$.
It is easy to see that the quadruple 
$$(V^{(0)}, \psi_2, \CR_1, \CC_2)$$
satisfies all the conditions in Lemma \ref{lem1}. 
Hence, applying Lemma \ref{lem1} to the quadruple
$(V^{(0)}, \psi_2, \CR_1, \CC_2)$,
the integral 
\begin{equation}\label{prop3equ6}
    \int_{[V^{(1)}]}
    \sum_{\xi_1 \in F^1} \omega_{\psi^{-1}}(\ell(\bar{p}_2,\bar{z})\wt{m})\varphi_{7}((0,\xi_1))
\end{equation}
\begin{equation*}
   \wt{E}_{1}(v^{(1)}\wt{m}w,\wt{\phi})
   \psi_{2}(v^{(1)})dv^{(1)}
\end{equation*}
is equal to
\begin{equation}\label{prop3equ7}
\int_{\CC_2(\BA)}
        \int_{[\CR_1]}
    \int_{[V^{(0)}]}
    \sum_{\xi_1 \in F^1} \omega_{\psi^{-1}}(\ell(\bar{p}_2,\bar{z})\ell(\bar{p}_1)\wt{m})\varphi_{7}((0,\xi_1))
\end{equation}
\begin{equation*}
   \wt{E}_{1}(v^{(0)}v^-(p_1)v\wt{m}w,\wt{\phi})\psi_{2}(v^{(0)})dv^{(0)}dp_1dv.
\end{equation*}
Hence, the integral \eqref{prop3equ5} is equal to 
\begin{equation}\label{prop3equ8}
\int_{[M]}\Phi^c(m,s,\phi)
\lvert \det a \rvert^{-\frac{21}{2}}
\int_{\CC_2(\BA)}
        \int_{[\CR_1]}
    \int_{[V^{(0)}]}
    \sum_{\xi_1 \in F^1} \omega_{\psi^{-1}}(\ell(\bar{p}_2,\bar{z})\ell(\bar{p}_1)\wt{m})
\end{equation}
\begin{equation*}
   \varphi_{7}((0,\xi_1))\wt{E}_{1}(v^{(0)}v^-(p_1)v\wt{m}w,\wt{\phi})\psi_{2}(v^{(0)})dv^{(0)}dp_1dvdm.
\end{equation*}

Since $[M]$ normalizes the group $\CC_2(\BA)$, after changing of variables, we obtain the following non-vanishing integral
\begin{equation}\label{prop3equ9}
\int_{\CC_2(\BA)}
\int_{[M]}\Phi^c(m,s,\phi)
\lvert \det a \rvert^{-\frac{23}{2}}
        \int_{[\CR_1]}
    \int_{[V^{(0)}]}
    \sum_{\xi_1 \in F^1} \omega_{\psi^{-1}}(\ell(\bar{p}_2,\bar{z})\ell(\bar{p}_1)\wt{m})
\end{equation}
\begin{equation*}
   \varphi_{7}((0,\xi_1))\wt{E}_{1}(v^{(0)}v^-(p_1)\wt{m}w (w^{-1}vw),\wt{\phi})\psi_{2}(v^{(0)})dv^{(0)}dp_1dmdv.
\end{equation*}
Therefore, as an inner integral, the following integral is non-vanishing
\begin{equation}\label{prop3equ10}
\int_{[M]}\Phi^c(m,s,\phi)
\lvert \det a \rvert^{-\frac{23}{2}}
    \int_{[\CR_1]}
    \int_{[V^{(0)}]}
    \sum_{\xi_1 \in F^1} \omega_{\psi^{-1}}(\ell(\bar{p}_2,\bar{z})\ell(\bar{p}_1)\wt{m})
\end{equation}
\begin{equation*}
   \varphi_{7}((0,\xi_1))\wt{E}_{1}(v^{(0)}v^-(p_1)\wt{m}wg,\wt{\phi})\psi_{2}(v^{(0)})dv^{(0)}dp_1dm,
\end{equation*}
where $g = w^{-1}vw$, for some $v \in \CC_2(\BA)$.

Now, we consider the multiple integral over $[\CR_1]$ and $[V^{(0)}]$, and 
apply Lemma \ref{lem1} to exchange the roots  
from 
$\CR_1$ to $\CC_2$.
Precisely, applying Lemma \ref{lem1} to the quadruple
$(V^{(0)}, \psi_2, \CC_2, \CR_1)$ (which also satisfies all the conditions in Lemma \ref{lem1}), 
the integral 
\begin{equation}\label{prop3equ11}
        \int_{[\CR_1]}
    \int_{[V^{(0)}]}
    \sum_{\xi_1 \in F^1} \omega_{\psi^{-1}}(\ell(\bar{p}_2,\bar{z})\ell(\bar{p}_1)\wt{m})\varphi_{7}((0,\xi_1))
\end{equation}
\begin{equation*}
   \wt{E}_{1}(v^{(0)}v^-(p_1)\wt{m}wg,\wt{\phi})\psi_{2}(v^{(0)})dv^{(0)}dp_1
\end{equation*}
is equal to
\begin{equation}\label{prop3equ12}
\int_{\CR_1(\BA)}
    \int_{[V^{(1)}]}
    \sum_{\xi_1 \in F^1} \omega_{\psi^{-1}}(\ell(\bar{p}_2,\bar{z})\ell(\bar{p}_1)\wt{m})\varphi_{7}((0,\xi_1))
\end{equation}
\begin{equation*}
   \wt{E}_{1}(v^{(1)}v^-(p_1)\wt{m}wg,
   \wt{\phi})
   \psi_{2}(v^{(1)})dv^{(1)}dp_1.
\end{equation*}

Hence, the integral \eqref{prop3equ10} is equal to 
\begin{equation}\label{prop3equ15}
\int_{[M]}\Phi^c(m,s,\phi)
\lvert \det a \rvert^{-\frac{23}{2}}
    \int_{\CR_1(\BA)}
    \int_{[V^{(1)}]}
    \sum_{\xi_1 \in F^1} \omega_{\psi^{-1}}(\ell(\bar{p}_2,\bar{z})\ell(\bar{p}_1)\wt{m})\varphi_{7}((0,\xi_1))
\end{equation}
\begin{equation*}
   \wt{E}_{1}(v^{(1)}v^-(p_1)\wt{m}wg,\wt{\phi})
   \psi_{2}(v^{(1)})dv^{(1)}dp_1dm.
\end{equation*}
Since $[M]$ normalizes $\CR_1(\BA)$, 
after changing of variables, we obtain the following non-vanishing integral
\begin{equation}\label{prop3equ18}
    \int_{\CR_1(\BA)}
\int_{[M]}\Phi^c(m,s,\phi)
\lvert \det a \rvert^{-\frac{21}{2}}
    \int_{[V^{(1)}]}
    \sum_{\xi_1 \in F^1} \omega_{\psi^{-1}}(\ell(\bar{p}_2,\bar{z})\wt{m}\ell(\bar{p}_1))
\end{equation}
\begin{equation*}
   \varphi_{7}((0,\xi_1))\wt{E}_{1}(v^{(1)}\wt{m}v^-(p_1)wg,\wt{\phi})\psi_{2}(v^{(1)})dv^{(1)}dmdp_1.
\end{equation*}
By choosing appropriate $\varphi_2 \in \CS(\BA^{6})$, the integral \eqref{prop3equ18} is non-vanishing if and only if the following integral is non-vanishing
\begin{equation}\label{prop3equ17}
    \int_{\CR_1\CR_2(\BA)}
\int_{[M]}\Phi^c(m,s,\phi)
\lvert \det a \rvert^{-\frac{21}{2}}
    \int_{[V^{(1)}]}
    \sum_{\xi_1 \in F^1} \omega_{\psi^{-1}}(\ell(\bar{p}_2,\bar{z})\wt{m}\ell(\bar{p}_1))
\end{equation}
\begin{equation*}
   \varphi_{7}((0,\xi_1))\wt{E}_{1}(v^{(1)}\wt{m}v^-(p_1)wg,\wt{\phi})\psi_{2}(v^{(1)})dv^{(1)}dmdp_1.
\end{equation*}

After taking Fourier expansion of $\wt{E}_1$ along $[\CC_1]$, arguing as in \cite[Section 4]{GJR04}, the integral \eqref{prop3equ17}
is equal to 
\begin{equation}\label{prop3equ14}
\int_{\CR_1\CR_2(\BA)}
\int_{[M]}\Phi^c(m,s,\phi)
\lvert \det a \rvert^{-\frac{21}{2}}
    \int_{[V_{3,1}]}
    \sum_{\xi_1 \in F^1} \omega_{\psi^{-1}}(\ell(\bar{p}_2,\bar{z})\wt{m}\ell(\bar{p}_1))
\end{equation}
\begin{equation*}
   \varphi_{7}((0,\xi_1))\wt{E}_{1,P_{6,3}}(v\wt{m}v^-(p_1)wg,\wt{\phi})\psi_{2}(v)dvdmdp_1.
\end{equation*}
Then, following the calculations from (4.29) -- (4.35) in \cite{GJR04}, we obtain that the following integral 
\begin{equation*}
  \frak{c} \int_{\Mat_{2,6}(\BA)\times\PM_{6,1}}
\phi_{\pi_\psi(\tiltau)\otimes\sig}(m)
{\CF}^\psi(R(g)\tilphi_{\pi_\psi(\tiltau)\otimes\tiltau})
(mv^-(p_1)w)dmdp_1
\end{equation*}
does not vanish, where $R(g)$ is the right translation operator. Hence, for  $\phi_{\pi_\psi(\tiltau)\otimes\sig} \in \CA^{ss}_{P_{2r,l}, \pi_{\psi}(\wt{\tau}) \otimes \sigma}$,  $R(g)\wt{\phi}_{\pi_\psi(\tiltau)\otimes\wt{\tau}} \in \wt{\CA}_{\wt{P}_{2r,r}, \pi_{\psi}(\wt{\tau}) \otimes \wt{\tau}, 0}$ would be suffice, in order for 
the integral \eqref{sec0equ4}
to be non-vanishing.

This completes the proof of Proposition \ref{prop2} for the special case $l=1$, $r=3$. \qed

\subsection{Proof of Proposition \ref{prop2}: the general case}\label{general case}

In this section, we prove Proposition \ref{prop2} for the general case. 

Again, we start from the non-vanishing integral \eqref{sec0equ20}. 
Reversing the calculations from (4.29) -- (4.35) as in \cite[Section 4]{GJR04}, the integral \eqref{sec0equ20} is equal to 
$\frac{1}{\frak{c}}$
times the residue at $s=\frac{1}{2}$ of 
\begin{equation}\label{sec0equ8}
\int_{[M]}\Phi^c(m,s,\phi) 
\lvert \det a \rvert^{-3r-l-\frac{1}{2}}
    \int_{[V_{r,l}]}
    \sum_{\xi_l \in F^l} \omega_{\psi^{-1}}(\ell(\bar{p}_2,\bar{z})\wt{m})\varphi_{2r+l}((0,\xi_l))
\end{equation}
\begin{equation*}
   \wt{E}_{1,\wt{P}_{2r,r}}(v\wt{m}w,\wt{\phi})\psi_{r-l}(v)dvdm,
\end{equation*}
where again $\frak{c}$
is as in \eqref{sec2equ5}. 

Reversing the step of taking Fourier expansion of $\wt{E}_1$ along $[\CC_1]$, the integral \eqref{sec0equ8} is equal to 
\begin{equation}\label{sec0equ9}
\int_{[M]}\Phi^c(m,s,\phi)
\lvert \det a \rvert^{-3r-l-\frac{1}{2}}
    \int_{[V^{(1)}]}
    \sum_{\xi_l \in F^l} \omega_{\psi^{-1}}(\ell(\bar{p}_2,\bar{z})\wt{m})\varphi_{2r+l}((0,\xi_l))
\end{equation}
\begin{equation*}
   \wt{E}_{1}(v^{(1)}\wt{m}w,\wt{\phi})\psi_{r-l}(v^{(1)})dv^{(1)}dm.
\end{equation*}
Recall that $V^{(1)}$ consists of elements of the type
$$v^{(1)}=\begin{pmatrix}
I_{2r} & q & y & p_3^* & z'\\
    & n & p_2 & z & p_3\\
    &   & I_{2l} & p_2^* & y^*\\
    &   &     & n^* & q^*\\
    &   &     &     & I_{2r}
\end{pmatrix},$$
where $q \in \Mat_{2r,r-l}$ with the first column being zero. 

Recall that $V^{(0)}$ consists of elements of the type
$$v^{(0)}=\begin{pmatrix}
I_{2r} & 0 & y & p_3^* & z'\\
    & n & p_2 & z & p_3\\
    &   & I_{2l} & p_2^* & y^*\\
    &   &     & n^* & 0\\
    &   &     &     & I_{2r}
\end{pmatrix}.$$
And for $1 \leq t \leq r-l$,
$$\CC_t = \left\{
\begin{pmatrix}
I_{2r} & q & 0 & 0 & 0\\
    & I_{r-l} & 0 & 0 & 0\\
    &   & I_{2l} & 0 & 0\\
    &   &     & I_{r-l} & q^*\\
    &   &     &     & I_{2r}
\end{pmatrix}: q \in \Mat_{2r,r-l}, q_{i,j}=0, j \neq t \right\},$$
$$\CR_t = \left\{
\begin{pmatrix}
I_{2r} & 0 & 0 & 0 & 0\\
  p  & I_{r-l} & 0 & 0 & 0\\
    &   & I_{2l} & 0 & 0\\
    &   &     & I_{r-l} & 0\\
    &   &     &   p^*  & I_{2r}
\end{pmatrix}: p \in \Mat_{r-l,2r}, q_{i,j}=0, i \neq t \right\}.$$

Note that $V^{(1)}=\prod_{i=2}^{r-l}\CC_iV^{(0)}$. Next, we consider the integral over $[V^{(1)}]$ and 
apply Lemma \ref{lem1} repeatedly to exchange roots from $\prod_{i=2}^{r-l}\CC_i$ to 
$\prod_{j=1}^{r-l-1}\CR_j$.
First, one can see that the quadruples 
$$\left(\prod_{j=1}^{t-1}\CR_j\prod_{i=t+2}^{r-l}\CC_i V^{(0)}, \psi_{r-l}, \CR_t, \CC_{t+1}\right), 1 \leq t \leq r-l-1,$$
satisfy all the conditions in Lemma \ref{lem1}.
Hence, applying Lemma \ref{lem1} to the
following ordered sequence of quadruples
$$\left(\prod_{i=3}^{r-l}\CC_i V^{(0)}, \psi_{r-l}, \CR_1, \CC_{2}\right),$$
$$\left(\CR_1\prod_{i=4}^{r-l}\CC_i V^{(0)}, \psi_{r-l}, \CR_2, \CC_{3}\right),$$
$$\cdots$$
$$\left(\prod_{j=1}^{t-1}\CR_j\prod_{i=t+2}^{r-l}\CC_i V^{(0)}, \psi_{r-l}, \CR_t, \CC_{t+1}\right),$$
$$\cdots$$
$$\left(\prod_{j=1}^{r-l-2}\CR_jV^{(0)}, \psi_{r-l}, \CR_{r-l-1}, \CC_{r-l}\right),$$
the integral 
\begin{equation}\label{sec0equ10}
    \int_{[V^{(1)}]}
    \sum_{\xi_l \in F^l} \omega_{\psi^{-1}}(\ell(\bar{p}_2,\bar{z})\wt{m})\varphi_{2r+l}((0,\xi_l))
\end{equation}
\begin{equation*}
   \wt{E}_{1}(v^{(1)}\wt{m}w,\wt{\phi})\psi_{r-l}(v^{(1)})dv^{(1)}
\end{equation*}
is equal to 
\begin{equation}\label{sec0equ11}
\int_{\prod_{i=2}^{r-l}\CC_i(\BA)}
\int_{[\prod_{j=1}^{r-l-1}\CR_j]}
    \int_{[V^{(0)}]}
    \sum_{\xi_l \in F^l} \omega_{\psi^{-1}}(\ell(\bar{p}_2,\bar{z})\ell(\bar{p}_1)\wt{m})\varphi_{2r+l}((0,\xi_l))
\end{equation}
\begin{equation*}
   \wt{E}_{1}(v^{(0)}v^-(p_1)v\wt{m}w,\wt{\phi})\psi_{r-l}(v^{(0)})dv^{(0)}dp_1dv.
\end{equation*}

Since $M$ normalizes the group $\prod_{i=2}^{r-l}\CC_i$, after changing of variables, we obtain the following non-vanishing integral
\begin{equation}\label{sec0equ12}
\int_{\prod_{i=2}^{r-l}\CC_i(\BA)}
\int_{[M]}\Phi^c(m,s,\phi)
\lvert \det a \rvert^{-4r+\frac{1}{2}}
\int_{[\prod_{j=1}^{r-l-1}\CR_j]}
    \int_{[V^{(0)}]}
    \sum_{\xi_l \in F^l} \omega_{\psi^{-1}}(\ell(\bar{p}_2,\bar{z})\ell(\bar{p}_1)\wt{m})
\end{equation}
\begin{equation*}
   \varphi_{2r+l}((0,\xi_l))\wt{E}_{1}(v^{(0)}v^-(p_1)\wt{m}w (w^{-1}vw),\wt{\phi})\psi_{r-l}(v^{(0)})dv^{(0)}dp_1dmdv.
\end{equation*}
Therefore, as an inner integral, the following integral is non-vanishing
\begin{equation}\label{sec0equ13}
\int_{[M]}\Phi^c(m,s,\phi)
\lvert \det a \rvert^{-4r+\frac{1}{2}}
\int_{[\prod_{j=1}^{r-l-1}\CR_j]}
    \int_{[V^{(0)}]}
    \sum_{\xi_l \in F^l} \omega_{\psi^{-1}}(\ell(\bar{p}_2,\bar{z})\ell(\bar{p}_1)\wt{m})
\end{equation}
\begin{equation*}
   \varphi_{2r+l}((0,\xi_l))\wt{E}_{1}(v^{(0)}v^-(p_1)\wt{m}wg,\wt{\phi})\psi_{r-l}(v^{(0)})dv^{(0)}dp_1dm,
\end{equation*}
where $g = w^{-1}vw$, for some $v \in \prod_{i=2}^{r-l}\CC_i(\BA)$.  

Now, we consider the multiple integral over $[\prod_{j=1}^{r-l-1}\CR_j]$ and $[V^{(0)}]$, and 
apply Lemma \ref{lem1} repeatedly to exchange roots 
from 
$\prod_{j=1}^{r-l-1}\CR_j$ to  $\prod_{i=2}^{r-l}\CC_i$. 
One can see that the quadruples 
$$\left(\prod_{i=2}^{t}\CC_{m}\prod_{j=t+1}^{r-l-1}\CR_j V^{(0)}, \psi_{r-l},  \CC_{t+1}, \CR_t\right), 1 \leq t \leq r-l-1,$$
also satisfy all the conditions in Lemma \ref{lem1}.
Hence, applying Lemma \ref{lem1} to the
following ordered sequence of quadruples
$$\left(\prod_{j=2}^{r-l-1}\CR_j V^{(0)}, \psi_{r-l}, \CC_{2}, \CR_1\right),$$
$$\left(\CC_2\prod_{j=3}^{r-l-1}\CR_j V^{(0)}, \psi_{r-l}, \CC_{3}, \CR_2\right),$$
$$\cdots$$
$$\left(\prod_{i=2}^{t}\CC_{i}\prod_{j=t+1}^{r-l-1}\CR_n V^{(0)}, \psi_{r-l},  \CC_{t+1}, \CR_t\right),$$
$$\cdots$$
$$\left(\prod_{i=2}^{r-l-1}\CC_{i}V^{(0)}, \psi_{r-l}, \CC_{r-l}, \CR_{r-l-1}\right),$$
the integral 
\begin{equation}\label{sec0equ14}
\int_{[\prod_{j=1}^{r-l-1}\CR_j]}
    \int_{[V^{(0)}]}
    \sum_{\xi_l \in F^l} \omega_{\psi^{-1}}(\ell(\bar{p}_2,\bar{z})\ell(\bar{p}_1)\wt{m})\varphi_{2r+l}((0,\xi_l))
\end{equation}
\begin{equation*}
   \wt{E}_{1}(v^{(0)}v^-(p_1)\wt{m}wg,\wt{\phi})\psi_{r-l}(v^{(0)})dv^{(0)}dp_1
\end{equation*}
is equal to 
\begin{equation}\label{sec0equ15}
\int_{\prod_{j=1}^{r-l-1}\CR_j(\BA)}
\int_{[\prod_{i=2}^{r-l}\CC_i]}
    \int_{[V^{(0)}]}
    \sum_{\xi_l \in F^l} \omega_{\psi^{-1}}(\ell(\bar{p}_2,\bar{z})\ell(\bar{p}_1)\wt{m})\varphi_{2r+l}((0,\xi_l))
\end{equation}
\begin{equation*}
   \wt{E}_{1}(v^{(0)}vv^-(p_1)\wt{m}wg,\wt{\phi})\psi_{r-l}(v^{(0)})dv^{(0)}dvdp_1.
\end{equation*}

Hence, the integral \eqref{sec0equ13} becomes
\begin{equation}\label{sec0equ16}
\int_{[M]}\Phi^c(m,s,\phi)
\lvert \det a \rvert^{-4r+\frac{1}{2}}
\int_{\prod_{j=1}^{r-l-1}\CR_j(\BA)}
    \int_{[V^{(1)}]}
    \sum_{\xi_l \in F^l} \omega_{\psi^{-1}}(\ell(\bar{p}_2,\bar{z})\ell(\bar{p}_1)\wt{m})
\end{equation}
\begin{equation*}
   \varphi_{2r+l}((0,\xi_l))\wt{E}_{1}(v^{(1)}v^-(p_1)\wt{m}wg,\wt{\phi})\psi_{r-l}(v^{(1)})dv^{(1)}dp_1dm.
\end{equation*}
Since $M$ normalizes the group $\prod_{j=1}^{r-l-1}\CR_j$, after changing of variables, we obtain the following non-vanishing integral
\begin{equation}\label{sec0equ17}
\int_{\prod_{j=1}^{r-l-1}\CR_j(\BA)}
\int_{[M]}\Phi^c(m,s,\phi)
\lvert \det a \rvert^{-3r-l-\frac{1}{2}}
    \int_{[V^{(1)}]}
    \sum_{\xi_l \in F^l} \omega_{\psi^{-1}}(\ell(\bar{p}_2,\bar{z})\wt{m}\ell(\bar{p}_1))
\end{equation}
\begin{equation*}
   \varphi_{2r+l}((0,\xi_l))\wt{E}_{1}(v^{(1)}\wt{m}v^-(p_1)wg,\wt{\phi})\psi_{r-l}(v^{(1)})dv^{(1)}dmdp_1.
\end{equation*}
By choosing appropriate $\varphi_2 \in \CS(\BA^{2r})$, the following integral is also non-vanishing
\begin{equation}\label{sec0equ19}
\int_{\prod_{j=1}^{r-l}\CR_j(\BA)}
\int_{[M]}\Phi^c(m,s,\phi)
\lvert \det a \rvert^{-3r-l-\frac{1}{2}}
    \int_{[V^{(1)}]}
    \sum_{\xi_l \in F^l} \omega_{\psi^{-1}}(\ell(\bar{p}_2,\bar{z})\wt{m}\ell(\bar{p}_1))
\end{equation}
\begin{equation*}
   \varphi_{2r+l}((0,\xi_l))\wt{E}_{1}(v^{(1)}\wt{m}v^-(p_1)wg,\wt{\phi})\psi_{r-l}(v^{(1)})dv^{(1)}dmdp_1.
\end{equation*}

After taking Fourier expansion of $\wt{E}_1$ along $[\CC_1]$, arguing as in \cite[Section 4]{GJR04}, the integral \eqref{sec0equ19} is equal to 
\begin{equation}\label{sec0equ18}
\int_{\prod_{j=1}^{r-l}\CR_j(\BA)}
\int_{[M]}\Phi^c(m,s,\phi)
\lvert \det a \rvert^{-3r-l-\frac{1}{2}}
    \int_{[V_{r,l}]}
    \sum_{\xi_l \in F^l} \omega_{\psi^{-1}}(\ell(\bar{p}_2,\bar{z})\wt{m}\ell(\bar{p}_1))
\end{equation}
\begin{equation*}
   \varphi_{2r+l}((0,\xi_l))\wt{E}_{1, \wt{P}_{2r,r}}(v\wt{m}v^-(p_1)wg,\wt{\phi})\psi_{r-l}(v)dvdmdp_1.
\end{equation*}
Then, following the calculations from (4.29) -- (4.35) as in \cite[Section 4]{GJR04}, we obtain that the following integral
$$
\frak{c} \int_{\Mat_{r-l,2r}(\BA)\times\PM_{2r,l}}
\phi_{\pi_\psi(\tiltau)\otimes\sig}(m)
{\CF}^\psi(R(g)\tilphi_{\pi_\psi(\tiltau)\otimes\tiltau})
(mv^-(p_1)w)dmdp_1,
$$
does not vanish, where $R(g)$ is the right translation operator. Hence, for  $\phi_{\pi_\psi(\tiltau)\otimes\sig} \in \CA^{ss}_{P_{2r,l}, \pi_{\psi}(\wt{\tau}) \otimes \sigma}$,  $R(g)\wt{\phi}_{\pi_\psi(\tiltau)\otimes\wt{\tau}} \in \wt{\CA}_{\wt{P}_{2r,r}, \pi_{\psi}(\wt{\tau}) \otimes \wt{\tau}, 0}$ would be suffice, in order for 
the integral \eqref{sec0equ4}
to be non-vanishing.

This completes the proof of Proposition \ref{prop2}. \qed

\textbf{Acknowledgement.}
We thank A. Ichino, E. Lapid, D. Soudry, S. Yamana, and L. Zhang for helpful 
discussions on the issue related to the previous proof of Proposition 5.3 in  \cite{GJR04}.

%------------------------------------------------------
\end{document}